\documentstyle[12pt]{article}
\begin{document}
\title{DIFFERENTIAL CALCULUS \\
ON A THREE-PARAMETER OSCILLATOR ALGEBRA}
\author{Mich\`{e}le IRAC-ASTAUD \thanks{e-mail : mici@ccr.jussieu.fr}\\ 
Laboratoire de Physique Th\'{e}orique et Math\'{e}matique\\
Universit\'{e} Paris VII\\2 place Jussieu F-75251 Paris Cedex 05, FRANCE}
\date{}
\maketitle
to be published in Reviews in Mathematical Physics
\begin{abstract}
Two differential calculi are developped on an algebra generalizing 
the usual q-oscillator algebra and involving three generators 
and three parameters. They are shown to be invariant under 
the same quantum group that is extended to a ten-generator Hopf algebra. 
We discuss the special case where it reduces to a deformation 
of the invariance group of the Weyl-Heisenberg algebra for which 
we prove the existence of a constraint between the values of the parameters.
\end{abstract}
\section{Introduction}

Generalizing the differential geometry on Lie groups and manifolds, 
the differential calculus on quantum groups and quantum spaces 
was developed in many interesting papers , (see for example 
\cite {Woro},\cite{Ber},\cite{wz},\cite{zumino}). An exterior 
differential,  one-forms and partial 
derivatives are defined  on a quadratic algebra, 
that is, a free associative algebra generated 
by variables satisfying quadratic 
commutation relations. In the present paper, 
as in \cite{nous} and \cite{nous1}, we consider 
the physically important case where   
the quadratic algebra  is a  deformation 
of the Weyl-Heisenberg algebra, in order to obtain 
a differential calculus on the algebra of the 
observables of a quantum system. 

The Weyl-Heisenberg algebra  can be seen as 
 a free associative algebra generated by three  
 variables $x^i$ satisfying the  quadratic relations :
\begin{equation}
R : \, \left\{
\begin{array}{rcl}
x^1 x^2 - x^2 x^1 - s(x^3)^2 &=&0\\
x^1 x^3 - x^3 x^1 &=& 0  \\
x^2 x^3 - x^3 x^2 &=& 0 \\
\end{array}
\right.
\label{1}
\end{equation}
where $\hat{a}= x^1, \quad 
\hat{a}^{\dagger}=  x^2$ and $x^3$ is the identity.
This algebra is  denoted by $C<x>/R$.

This space is invariant under the seven-dimensional Lie  subgroup of 
 $GL(3)$, $G$, constituted by the matrices 
 $T$ such as $t^3_1 =t^3_2 =0$ and $t^1_1t^2_2-t^1_2t^2_1=(t^3_3)^2$.

The quadratic deformation of the Weyl-Heisenberg 
algebra considered in this paper, as in  \cite{nous} 
and \cite{nous1}, is an algebra,  generalizing $C<x>/R$, 
obtained by replacing the
  relations (\ref{1}) by
\begin{equation}
R_{xx} : \, 
\left\{
\begin{array}{rcl}
x^1 x^2 - q \, x^2 x^1 -s (x^3)^2 & = & 0 \\
x^1x^3- u \, x^3 x^1 &  = & 0 \\
x^2x^3-u^{-1}\, x^3 x^2 & = & 0
\end{array}
\right.
\label{xx}
\end{equation}
When $q=u^{-2}$,  putting
$x^1=\hat{a}, \quad 
x^2= \hat{a}^{\dagger}$ and $x^3=q^{-\frac{N}{2}}$, 
the algebra $C<x>/R_{xx}$ is identified to the 
q-oscillator algebra \cite{macfarlane}\cite{biedenharn}. 
When $s=0$, $C<x>/R_{xx}$ is the three-dimensional 
quantum plane \cite{manin}. 

In  \cite{nous} and \cite{nous1}, we establish that a 
differential calculus on $C<x>/R_{xx}$  invariant under a seven-dimensional 
quantum group, deformation of $G$, only exists if $q= u^2$ and then is  
unique. As the  constraint, $q= u^2$,
 eliminates the important case when the quantum space is 
 the q-oscillator algebra, 
the aim of this paper is to obtain
   a   differential calculus on $C<x>/R_{xx}$  not restricted, 
   a priori, by a requirement of invariance and valid for arbitrary 
   values of the deformation parameters $q,u$ and $s$.

 In section 2, we  determine  two different sets of  consistent 
 quadratic relations  between variables, differentials and derivatives. 
 In section 3, we prove that all these relations are invariant under 
 the same quantum group that is extended to a ten-generator Hopf-algebra.   
 When we assume that the variables $x^1$ and $x^2$  are  mutually adjoint 
 and that $x^3$ is self-adjoint, the quantum group is endowed with a 
 structure of Hopf-star-algebra. Finally we discuss 
the particular case where the quantum group is a deformation of  $G$ and 
we recover the constraint on the parameters \cite{nous}, \cite{nous1}. 
In this particular case, the  ten-generator Hopf-star-algebra  contains a  
eight-generator subalgebra , both algebras leave invariant the commutation 
relations defining the differential calculus.

\section{Differential Calculus} 
Following the usual method \cite{wz}\cite{zumino}, we add  to the 
free algebra $C<x>$,
 three generators $\xi^i,\, i=1,2,3$, identified to the one-forms.  
 We define the exterior 
differential operator $d $   in $C<x,\xi>$ such as :
 
\noindent $ d(x^i) = \xi^i$,
  
\noindent $d$ is linear,

\noindent $ d^2=0$

\noindent and  
$ d$  satisfies the graded  Leibniz rule :
\begin{equation}
d(fg) = (df)g + (-1)^k f(dg)
\end{equation} 
where $ f,g \in {C<x,\xi>}$ and $f$ is of degree $k$. 
Then,  the 
partial derivatives $\partial_i$ are defined by~:

$$ d \equiv \xi^{i} \; \partial_i.$$

The commutation relations $R_{x\xi}$ between the variables and the 
differentials and those between the partial derivatives $\partial_k$ 
and the forms  $\xi^l$ are assumed to be quadratic \cite{wz}\cite{zumino} :
\begin {equation} 
R_{x\xi} \quad: \quad  x^k\xi^l=C^{kl}_{mn} \xi^m x^n 
\label{31}
\end{equation}
and
\begin{equation}
R_{\partial \xi} \quad : \quad
\partial_k \xi^l=K^{lm}_{kn}\xi^n \partial_m
\label{xder} 
\end{equation}
 By applying operator $d$ to (\ref{31})
on the left, we get :
\begin{equation} 
R_{\xi\xi} \quad :\quad \xi^k \xi^l=-C^{kl}_{mn} \xi^m\xi^n ,\quad \! \! 
\label{32}
\end{equation} 
The above definition implies :
$\partial_l(x^k)=\delta^k_l$.
Applying the Leibniz rule on $x^k \, f$, $f\in C<x>$, and taking into 
account
relations (\ref{31}) , we obtain the commutation relations 
between $\partial_i$ and $x^k$ \cite{wz}:
 \begin{equation}
R_{x\partial} \quad :\quad
\partial_lx^k=\delta^k_l+C^{km}_{ln}x^n\partial_m
\label{partialx}
\end{equation}

The combinations of the relations $R_{\partial \xi}, R_{x\partial}$ 
and $R_{x\xi}$ lead to
\begin{equation}
K  \! \! = \! \!  C^{-1}
\end{equation}
and to a sufficient condition of consistency, the Yang-Baxter equation :
\begin{equation}
(C\otimes 1)(1 \otimes C)(C\otimes 1) = (1 \otimes C)(C\otimes 1)
(1 \otimes C)
\label{yb1}
\end{equation}
where 1 is the identity of $GL(3)$.

From this point, we simplify the matrix $C$ assuming that 
$C^{ij}_{lm}$ is zero if $(i,j)\ne (l, m)$ or $(m,l)$ except 
$C^{12}_{33}$ and $C^{21}_{33}$.
We multiply  the relations $R_{xx}$  on the left by $\partial_i$. 
Using $R_{x\partial}$ , we commute $\partial_i$ to the right, and obtain 
several relations between the elements $C^{ij}_{kl}$. In particular,
\begin{equation}
\begin{array}{llll}
C^{12}_{12}&=q C^{21}_{12}-1,& C^{12}_{21}&=q C^{21}_{21} +q \\
C^{13}_{13}& =u C^{31}_{13} -1,& C^{13}_{31}&= u C^{31}_{31}+u\\
C^{32}_{23}&=u C^{23}_{23} +u,& C^{32}_{32}&=u C^{23}_{32}-1 \\
C^{12}_{33}& = q C^{21}_{33} +s C^{33}_{33} +s&& 
\end {array}
\label{coeff1}
\end{equation}
and 
\begin{equation}
C^{12}_{12} C^{21}_{21} = 0 ,\quad
 C^{13}_{13} C^{31}_{31} = 0 ,\quad
C^{23}_{23} C^{32}_{32} = 0 
\label{coeff2}
\end{equation}
We substitute in the equation (\ref{yb1}) a matrix $C$ satisfying 
(\ref{coeff1}) . This gives $27 \times 27$ new relations between the 
coefficients $C^{ij}_{kl}$. Solving these relations together with 
(\ref{coeff2}), we find only two different solutions 
: $C$ is equal to $\Omega$ or to its inverse, with  $\Omega$ defined by :
\begin{equation}
\Omega =\left(\begin{array}{c c c c c c c c c }
q/u^{2}&0&0&0&0&0&0&0&0\\
0&0&0&q^2/u^2&0&0&0&0&qs/u^2\\
0&0&0&0&0&0&q/u&0&0\\
0&q^{-1}&0&q/u^2-1&0&0&0&0&-s/q\\
0&0&0&0&q/u^2&0&0&0&0\\
0&0&0&0&0&q/u^2-1&0&1/u&0\\
0&0&1/u&0&0&0&q/u^2-1&0&0\\
0&0&0&0&0&q/u&0&0&0\\
0&0&0&0&0&0&0&0&q/u^2\\
\end{array}\right)  
\label{C1}
\end{equation}  
The matrices $\Omega$ and $\Omega ^{-1}$ have the same  eigenspaces 
that correspond to the variable quantum space defined by $R_{xx}$ 
and to the one-form quantum space defined by :
\begin{equation} 
R_{\xi\xi} : \,
\left
\{\begin{array}{cccccc}
(\xi^1)^2 & = & 0,\quad & (\xi^2)^2 & = & 0,\quad \\
(\xi^3)^2 & = & 0 \quad &\xi^2 \xi^1& = &-  u^{2}/q^2 \, \xi^1\xi^2 \\  
\xi^1 \xi^3&=&-q/u \, \xi^3 \xi^1, \quad 
& \xi^2 \xi^3 &= & - u/q \, \xi^3 \xi^2 
\end{array}
\right.
\label{xi1}
\end{equation}
When $C$ is equal to $\Omega$ and to its inverse, the eigenspaces 
of the transpose matrix $(C^{-1})^t$ are the same. The six-dimensional 
eigenspace
 is identified to the derivative quantum space and is defined by :
\begin{equation} 
R_{\partial \partial} : \,
\left
\{\begin{array}{ll}
 \partial_1 \partial_2 = & u^{2}/q^2 \,\partial_2\partial_1,\\  
\partial_1 \partial_3=&u/q \, \partial_3 \partial_1,\\
\partial_2\partial_3 = &  q/u \, \partial_3 \partial_2 
\end{array}
\right.
\label{partial1}
\end{equation}
The three-dimensional eigenspace corresponds to the covariant differential 
forms. We denote $R$ the set of relations (\ref{xx}), (\ref{xi1}) 
and (\ref{partial1}).
Corresponding to $C=\Omega$ or $C=\Omega^{-1}$, we obtain two  sets of 
relations $R_{x\xi}$, $R_{\partial \xi}$ and $R_{\partial x}$ and then 
two different differential calculi :

$\bullet$ The set of relations $R^{\Omega}$ associated to $\Omega$ is :

\noindent The commutation relations $R_{x\xi}^{\Omega}$ between the 
variables and the differentials~:
\begin{equation}
\begin{array}{llll}
 x^i \xi^i&= q/u^2\,\xi^i x^i,\quad i=1,2,3,&x^1 \xi^3&=q/ u \,\xi^3 x^1,\\
x^1 \xi^2&= q^2/u^2 \,\xi^2 x^1+ qs/u^2 \, \xi^3 x^3,&
x^3 \xi^2&= q/u \, \xi^2 x^3,\\
x^2 \xi^3&= (q/u^2-1) \, \xi^2 x^3+ 1/u \, 
\xi^3 x^2,&x^3 \xi^1&=( q/u^2-1) \, \xi^3 x^1+ 1/u \, \xi^1 x^3\\
x^2 \xi^1&= 1/q \, \xi^1 x^2+(q/u^2-1) \, \xi^2 x^1- &s/q \,& \xi^3 x^3,\\
\end{array}
\end{equation}
The commutation relations $R_{\partial \xi}^{\Omega}$ between the 
derivatives and the differentials~:
\begin{equation}
\begin{array}{llll}
\partial_3 \xi^3&= (u^2/q -1)\,\xi^2\partial_2+ u^2/q \, \xi^3\partial_2
&\partial_1 \xi^2&=  u^2/q^2 \,\xi^2\partial_1,\\
\partial_1 \xi^3&=  u/q \,\xi^3\partial_1,
&\partial_2 \xi^1&= q \, \xi^1\partial_2,\\
\partial_3 \xi^2&=  u/q \,\xi^2\partial_3- su^2/q^2 \,\xi^3\partial_1,
&\partial_2 \xi^3&= u \, \xi^3\partial_2,\\
\partial_3 \xi^1&= u \, \xi^1\partial_3+ s \, \xi^3\partial_2,
&\partial_2 \xi^2&=  u^2/q \, \xi^2\partial_2,\\
\partial_1 \xi^1&=  u^2/q \, \xi^1\partial_1+(u^2/q -1)\,
\ \xi^3\partial_3+&(u^2/q -1)&\,\ \xi^2\partial_2,
\end{array}
\end{equation}
The commutation relations $R_{x \partial}^{\Omega}$ between 
the derivatives and the variables~:
\begin{equation}
\begin{array}{llll}
  \partial_1 x^1&= 1+  q/u^2\, x^1\partial_1,&\partial_2 x^3&=  
  q/u \,x^3\partial_2,\\
\partial_3 x^3&= 1+  q/u^2\,x^3\partial_3+ (q/u^2-1) \, 
x^1\partial_1,&\partial_1 x^2&= 1/q \, x^2\partial_1,\\
\partial_3 x^1&=  q/u \,x^1\partial_3+  qs/u^2\, x^3\partial_2,
&\partial_2 x^1&=  q^2/u^2\,x^1\partial_2,\\
\partial_3 x^2&= 1/u \, x^2\partial_3-s/q \,x^3\partial_1,&\partial_1 
x^3&= 1/u \, x^3\partial_1,\\
\partial_2 x^2&= 1+ q/u^2\,x^2\partial_2+  (q/u^2-1) &\,x^1\partial_1+& 
(q/u^2-1) \, x^3\partial_3.\\
\end{array}
\end{equation}

$\bullet$ The set of relations associated to $\Omega^{-1}$ :

\noindent The commutation relations $R_{x\xi}^{\Omega^{-1}}$ between 
the variables and the differentials~:
\begin{equation}
\begin{array}{llll}
 x^i \xi^i&=  u^2/q \,\xi^i x^i,\quad i=1,2,3,&x^1 \xi^3&= ( u^2/q -1) \, 
 \xi^1x^3+ u \, \xi^3 x^1,\\
x^3 \xi^1&= u/q \, \xi^1x^3, &x^2\xi^1&= u^2/q^2 \, \xi^1x^2- su^2/q^2 
\,\xi^3x^3,\\
x^2 \xi^3&= u/q \, \xi^3 x^2 ,&x^3 \xi^2&= ( u^2/q -1) \,\xi^3 x^2+ u  
\,\xi^2 x^3,\\
x^1 \xi^2&= ( u^2/q -1) \,\xi^1 x^2 + q  \,\xi^2 x^1 + &s \,\xi^3x^3.&
\end{array}
\end{equation}
The commutation relations $R_{\partial \xi}^{\Omega^{-1}}$ between the 
derivatives and the differentials~:
\begin{equation}
\begin{array}{llll}
 \partial_1 \xi^1&= q/u^2\,\xi^1\partial_1,&\partial_3 \xi^2&= 1/u \, 
 \xi^2\partial_3- s/q \, \xi^3\partial_1,\\
\partial_1 \xi^3&= 1/u \, \xi^3\partial_1,
&\partial_2 \xi^1&= q^2/u^2\,\xi^1\partial_2,\\
\partial_2 \xi^3&= q/u \,\xi^3\partial_2,&\partial_3 \xi^1&= q/u \, 
\,\xi^1\partial_3+ sq/u^2\, \,\xi^3\partial_2,\\
\partial_1 \xi^2&= 1/q \, \xi^2\partial_1,
&\partial_3 \xi^3&= (q/u^2-1) \,\xi^1\partial_1+ q/u^2\, \xi^3\partial_3,\\
\partial_2 \xi^2&= (q/u^2-1) \,\xi^1\partial_1+&(q/u^2-1) \,
&\xi^3\partial_3+ q/u^2\,\xi^2\partial_2
\end{array}
\end{equation}
The commutation relations $R_{x \partial}^{\Omega^{-1}}$ between 
the derivatives and the variables~:
\begin{equation}
\begin{array}{llll}
\partial_2 x^2&=1+ u^2/q \,x^2\partial_2, 
&\partial_1 x^3&= u/q \,x^3\partial_1,\\
\partial_3 x^3&= 1+ u^2/q \,x^3\partial_3+ (u^2/q -1) \,x^2\partial_2,&
\partial_2 x^1&= q \, x^1\partial_2,\\
\partial_1 x^2&= u^2/q^2 \, x^2\partial_1,
&\partial_3 x^2&= u/q \,x^2\partial_3- su^2/q \,x^3\partial_1,\\
\partial_3 x^1&= u \,x^1\partial_3 +s \, x^3 \partial_2,
&\partial_2 x^3&= u \, x^3\partial_2,\\
\partial_1 x^1&= 1+ u^2/q \, x^1\partial_1+(u^2/q -1) 
\,x^2\partial_2+&(u^2/q -1) \,& x^3\partial_3.\\
\end{array}
\end{equation}
The two sets of relations $R_{xx},R_{\xi\xi},R_{\partial\partial}$ 
with $R^{\Omega}$ or $R^{\Omega^{-1}}$ define two quadratic algebras 
:$C<x,\xi,\partial>/R \cup R^{\Omega}$ and 
$C<x,\xi,\partial>/R \cup R^{\Omega^{-1}}$. In the following section, 
we investigate their invariance.

It is to be noted that all the  construction of the differential 
calculus is performed without using  the B-matrix associated with 
the variables \cite{wz} \cite{zumino}, and is the result solely  of 
the relations $R_{xx}, R_{x \xi}$ and $R_{\partial \xi}$ . Moreover, 
as a consequence of the construction, $B$ is found to be equal to $C$. 
\section{Quantum Group and Invariance}
The quantum matrix $T$ with nine non commuting elements defines a 
homomorphism on $C<x,\xi,\partial>$ \cite{rtf}. The  variables $x$ 
and  the differentials $\xi$ are transformed by $T$ and 
the derivatives $\partial$ are transformed by $(T^{-1})^ t$.

\noindent When the matrix $T$ satisfies 
\begin{equation}
{R}^{ji}_{kl} \, t^k_m t^l_n = t^j_l t^i_k \, {R}^{lk}_{mn}
\end{equation}
with $R=\Omega$ (resp. $R=\Omega^{-1}$), the  relations 
$R \cup R^{\Omega}$ (resp. $R \cup R^{\Omega^{-1}}$) are 
invariant, and therefore  this homomorphism maps $C<x,\xi,
\partial>/R \cup R^{\Omega}$ (resp. 
$C<x,\xi,\partial>/R \cup R^{\Omega^{-1}}$) in itself.
  It is easy to see that $R=\Omega$ and $R=\Omega^{-1}$ define 
  the same quantum matrix $T$, the elements of which satisfy the 
  following commutation relations 
$R_{tt} \quad :$ 
\begin{equation}
\begin{array}{lllll}
t^1_2t^1_1&=q^2/u^2\, t^1_1t^1_2,&t^2_2t^1_1&=t^1_1t^2_2-(u^2-q)/q^2\, 
t^1_2t^2_1- qs /u^2\, t^3_1t^3_2,\\
t^1_3t^1_2&= u /q \,t^1_2t^1_3,&t^2_1t^1_1&=1 /q \,t^1_1t^2_1 - s/q 
\,(t^3_1)^2,\\
t^1_3t^1_1&= q/u \,t^1_1t^1_3,&t^2_3t^1_1&= u/q \,t^1_1t^2_3-
(u^2-q)/q^2 \,t^1_3t^2_1- s/q \,t^3_3t^3_1,\\
t^3_2t^2_2&=u \,t^2_2t^3_2,&t^3_3t^1_1&=t^1_1t^3_3- (u^2-q)/(uq)\, 
t^1_3t^3_1,\\
t^3_1t^1_1&= 1/u\, t^1_1t^3_1,&t^3_2t^1_1&= q/u \,t^1_1t^3_2- 
(u^2-q)/u\, t^1_2t^3_1,\\
t^3_3t^1_2&=1/q \,t^1_2t^3_3,&t^2_2t^1_2&= 1/q \,t^1_2t^2_2- s/q 
\,(t^3_2)^2,\\
t^2_3t^2_2&= u/q \,t^2_2t^2_3,&t^2_3t^1_2&= u/q^2\, 
t^1_2t^2_3- s/q \,t^3_3t^3_2,\\
t^3_1t^1_2&= u/q^2\, t^1_2t^3_1,&t^2_1t^1_2&= u^2/q^3 
\,t^1_2t^2_1- s/q\, t^3_1t^3_2,\\
t^3_2t^2_1&= q^2/u \,t^2_1t^3_2,&t^3_2t^1_3&=t^1_3t^3_2- 
(u^2-q)/(uq)\, t^1_2t^3_3,\\
t^3_2t^1_2&=1/u \,t^1_2t^3_2, &t^2_3t^1_3&= 1/q\, t^1_3t^2_3- 
s/q \,(t^3_3)^2+ s/q \,t^1_1t^2_2- su^2/q^3\, t^1_2t^2_1,\\
t^3_1t^1_3&= 1/q\, t^1_3t^3_1,&t^2_2t^1_3&= u/q\,t^1_3t^2_2- 
(u^2-q)/q^2\, t^1_2t^2_3- s/u \,t^3_3t^3_2,\\
t^2_1t^2_2&= u^2/q^2\, t^2_2t^2_1,&t^3_3t^2_2&=t^2_2t^3_3+ 
(u^2-q)/u \,t^2_3t^3_2\\
t^2_3t^2_1&= q/u \,t^2_1t^2_3,&t^3_3t^2_3&=u \,t^2_3t^3_3+ 
sq/u \,t^2_1t^3_2-su\,t^2_2t^3_1,\\
t^3_2t^3_1&= q^2/u^2\, t^3_1t^3_2,&t^3_1t^2_2&= u/q \,t^2_2t^3_1 + 
(u^2-q)/u \,t^2_1t^3_2,\\
t^3_2t^3_3&= q/u \,t^3_3t^3_2,&t^3_1t^2_3&=t^2_3t^3_1+ (u^2-q)/u 
\,t^2_1t^3_3,\\
t^3_1t^2_1&=u \,t^2_1t^3_1,& t^2_1t^1_3&= u/q^2 \,t^1_3t^2_1 - 
su/q^2 \,t^3_3t^3_1\\
t^3_1t^3_3&= u/q \,t^3_3t^3_1,&t^3_3t^1_3&=1/u \,t^1_3t^3_3 +s/u\, 
t^1_1t^3_2-su/q^2 \,t^1_2t^3_1,\\
t^3_3t^2_1&=q \,t^2_1t^3_3,&t^3_2t^2_3&=q \,t^2_3t^3_2\\
\label{RT}
\end{array}
\end{equation}
Any elements  of $C<t>/R_{tt}$ can be written as a sum of ordered 
monomials $(t^1_1)^{k_1} (t^1_2)^{k_2}(t^1_3)^{k_3}
(t^2_2)^{k_4}(t^2_1)^{k_5}(t^2_3)^{k_6}(t^3_3)^{k_7}
(t^3_1)^{k_8}(t^3_2)^{k_9}$ by using the relations $R_{tt}$.
The inverse $T^{-1}$, of $T$ is equal to :
\begin{equation}
\left(\begin{array}{c c c}
t^2_2t^3_3-q/u \,t^2_3t^3_2&-q^2/u^2\, t^1_2t^3_3 + q^3/u^3\, 
t^1_3t^3_2&t^1_2t^2_3-q/u \,t^1_3t^2_2\\
-u^2/q^2\, t^2_1t^3_3 +u^3/q^3 \,t^2_3t^3_1&t^1_1t^3_3-u/q 
\,t^1_3t^3_1&-u^2/q^2 \,t^1_1t^2_3+ u^3/q^3 \,t^1_3t^2_1\\
t^2_1t^3_2 -u^2/q^2 \,t^2_2t^3_1&-q^2/u^2 \,t^1_1t^3_2 
+t^1_2t^3_1&t^1_1t^2_2 -u^2/q^2\, t^1_2t^2_1\\
\end{array}\right) D^{-1} 
\label{t1}
\end{equation}
with the determinant of $T$ equal to 
$$D =t^1_1t^2_2t^3_3+t^1_3t^2_1t^3_2+u^3/q^3 \,t^1_2t^2_3t^3_1- 
q/u \,t^1_1t^2_3t^3_2-u^2/q^2 \,t^1_2t^2_1t^3_3 -u^2/q^2\, 
t^1_3t^2_2t^3_1.$$ 
We can calculate and verify that $D$ is not a central 
element of $C<t>$ and therefore we have to add the generator 
$D^{-1}$ to this algebra. The commutation relations 
$R_{t D^{-1}} $ of $D^{-1}$ with all the $t^i_j$ are 
deduced from those of $D$ :
\begin{equation}
\begin{array}{llllll}
t^1_1 D^{-1}&=D^{-1} t^1_1,& t^1_2 D^{-1}&=u^2/q^4\, 
D^{-1}t^1_2,& t^1_3 D^{-1}&= u/q^2 \,D^{-1}t^1_3,\\
t^2_2 D^{-1}&=D^{-1} t^2_2,& t^2_1 D^{-1}&= q^2 \,
D^{-1}t^2_1,& t^2_3 D^{-1}&= u/q^2 \,D^{-1} t^2_3 ,\\
t^3_1 D^{-1} &= q^2/u \,D^{-1} t^3_1,& t^3_2 D^{-1}&=u/q^2\,  
D^{-1} t^3_2,& t^3_3 D^{-1}&= D^{-1}t^3_3.\\
\end{array}
\label{RD}
\end{equation}
The quotient algebra $C< t, D^{-1}>/R_{tt} \cup R_{t D^{-1}}$ 
is a Hopf algebra with the co-product $\Delta$, 
the co-unit $\epsilon$ and antipode $S$  defined by :  

\begin{eqnarray}
\Delta(T) \equiv T \otimes T, \quad 
\Delta(D^{-1}) \equiv D^{-1} \otimes D^{-1}
\label{26} \\
\epsilon (T,D^{-1}) \equiv (I,1), \quad 
S(T)\equiv T^{-1}, \quad   S(D) \equiv D^{-1}
\label{28}
\end{eqnarray}

When $x^1$ and $x^2$ are mutually adjoint and 
when $x^3$ is self-adjoint (for instance, 
in the case of the q-oscillator algebra ), 
the relations $R_{xx}$ are unchanged if 
the parameters are real. The action of $T$ respects 
this property if the quantum group is equipped with a 
star-operation that is an antihomomorphism such as 
\begin{equation}
\begin{array}{llllll}
&((t^i_j)^*)^*&=t^i_j,&\quad &\forall i,j&\\
(t^2_2)^*&=t^1_1,&\quad (t^2_1)^*&=t^1_2,&\quad
 (t^2_3)^*&=t^1_3,\\
&( t^3_1)^*&=t^3_2,&\quad (t^3_3)^*&=t^3_3&.
\end{array}
\label{star}
\end{equation}
These relations are consistent with (\ref{RT}) and (\ref{RD}) 
and the quantum group $C< t, D^{-1}>/R_{tt} \cup R_{t D^{-1}}$ 
acquires the structure of a Hopf-star-algebra.

Let us stress some properties of $C< t, D^{-1}>/R_{tt} \cup R_{t D^{-1}}$~:

$\bullet$ 
When $s=0$, the variable quantum space is the three-dimensional 
quantum plane, the resulting quantum group  corresponds to an original 
deformation of $GL(3)$ \cite{rtf} , \cite{schirr}.

$\bullet$ 
When $t^3_1$ and $t^3_2$ vanish, we obtain a deformation $G_{q s}$ 
of the subgroup $G$ of $GL(3)$.
      Two of the relations $R_{tt}$ give  : 
$$(u^2 -q)\, t^1_2 t^3_3=0,$$
 and 
 $$(u^2 -q)\, t^2_1 t^3_3=0.$$ 
implying that  $q$ is equal to $u^2$ if the algebra has no zero divisors. 
When $q=u^2$,
 the matrix $\Omega$ being equal to its inverse,  the two 
 differential calculi on $C<x>/R_{xx}$ reduce to one. This differential 
 calculus was previously obtained by a completely different method, 
 implying the uniqueness of the result \cite{nous1}.  All the 
 commutation relations are invariant under the ten-generator 
 quantum group $C< t, D^{-1}>/R_{tt} \cup R_{t D^{-1}}$ and under 
 a quantum group $G_{q s}$ deduced from the previous one  by putting 
 $t^3_1=t^3_2=0.$

$\bullet$
The relations (\ref{star}) are consistent with the condition 
$t^3_1=t^3_2=0$, and $G_{q s}$  is a Hopf-star-subalgebra 
of $C< t, D^{-1}>/R_{tt} \cup R_{t D^{-1}}$.
 
$\bullet$ 
In the case where $q=u^2$, we would point out that, if we add the 
generator $(t^3_3)^{-1}$ to the quantum group $G_{q s}$, the $T$-matrix 
can be written on the form $T^\prime \times t^3_3$ 
with $t^{\prime i}_j = t^i_j (t^3_3 )^{-1}$. All the elements 
$t^{\prime i}_j$ commute two by two. Nevertheless, they cannot 
be identified to C-numbers (and then $T^\prime$ to a matrix 
belonging to the initial subgroup $G$), because if they were 
C-numbers all the elements $t^i_j$ would be proportional to 
$t^3_3$ and this is impossible due to their commutation relations 
resulting from (\ref{RT}).
\section{Conclusion}
Two differential calculi can be associated with the three-parameter 
oscillator algebra $C<x>/R_{xx}$ and in particular with the q-oscillator 
algebra. They are invariant under the same quantum group 
$C< t, D^{-1}>/R_{tt} \cup R_{t D^{-1}}$ that is an original 
three-parameter deformation of $GL(3)$. When we assume that the 
variables $x^1$ and $x^2$  are  mutually adjoint and that $x^3$ 
is self adjoint, a star-operation is defined on the invariance 
quantum group that  then becomes a Hopf-star-algebra. Finally we 
consider the case where two generators, $t^3_1$ and $t^3_2$, are 
removed from the algebra $C< t, D^{-1}>/R_{tt} \cup R_{t D^{-1}}$ 
and we prove that the resulting differential calculus and quantum 
group of invariance do not exist for arbitrary values of the parameters, 
this is in agreement with a previous work \cite{nous}\cite{nous1}.

We would like to thank J.Bertrand for the many interesting 
discussions we have had.

 \end{document}